\documentclass[11pt]{amsart}

\usepackage{amsmath,amsthm}
\numberwithin{equation}{section}

\newtheorem{theorem}{Theorem}[section]

\theoremstyle{definition}

\theoremstyle{remark}
\newtheorem{remark}[theorem]{Remark}

\def\to{{\rightarrow}}
\def\R{\mathbb{R}}
\def\C{\mathbb{C}}
\def\e{\mbox{e}}

\title{A new algorithm for the volume of a convex polytope}

\thanks{Research partially supported the ECOS-Nord (France)-ANUIES (M\'exico)
Educational and Scientific Cooperation Program PM98M02}

\author{Jean B. Lasserre}
\address{LAAS-CNRS, 7 Avenue du Colonel Roche,
31077 Toulouse C\'{e}dex 4, France.}
\email{lasserre@laas.fr}

\author{Eduardo S. Zeron}
\address{Depto. Matem\'aticas, CIVESTAV-IPN,
Apdo.~Postal 14740, Mexico D.F. 07000, M\'exico.}
\email{eszeron@math.cinvestav.mx}

\date{}

\begin{document}

\begin{abstract}
We provide two algorithms 
for computing the volume of the convex polytope
$\Omega:=\{x\in \R^n_+ \,|\,Ax\leq b\}$, for $A\in\R^{m\times n}, b\in\R^n$.
Both algorithms have a $O(n^m)$  computational complexity which makes
them especially attractive for large $n$ and relatively small $m$, when
the other methods with $O(m^n)$ complexity fail.
The methodology which differs from previous existing methods uses
a Laplace transform technique that is well suited to the half-space
representation of $\Omega$.
\end{abstract}

\maketitle

\section{Introduction}

In this paper, we are interested in the exact computation
of the volume of the convex polytope
$\Omega :=\{x\in\R^n_+\,|\,Ax\leq b\}$, for some given matrix 
$A\in\R^{m\times n}$ and vector $b\in\R^m$.

Computing the volume of a convex polytope $\Omega$ is difficult. 
Basically, methods for exact computation of this volume can be 
classified according to whether one has a {\it half-space} 
representation of $\Omega$ as above, or a {\it vertex} representation, 
that is, when $\Omega$ is given by its list of vertices (triangulation
methods), or when both descriptions are available.
For instance, Lasserre 's algorithm \cite{lasserre} requires
a half-space description, whereas Delaunay 's triangulation (see
e.g. \cite{bueler}) or Von Hohenbalken's simplicial decomposition
\cite{von} require
the list of vertices. On the other hand, both Lawrence 's formula
\cite{lawrence} and Cohen and Hickey 's triangulation method \cite{cohen} require the double
(half-space and vertex) description of the polytope. 
For an updated review 
of the above methods and their
computational complexity, the interested reader is referred to
B\"ueler et al \cite{bueler}. In  particular, 
improved versions of some of the above algorithms are also described
in \cite{bueler}. The computational complexity is
also discussed in Dyer and Frieze \cite{dyer2}.
In a different spirit, Barvinok
\cite{barvinok} approximates the volume by computing the integral
of $\exp{\langle c ,x\rangle}$ over $\Omega$ for a small $c$. 
The latter integral reduces to evaluate at each vertex $v$ of 
$\Omega$, the integral of $\exp{\langle c ,x\rangle}$ over the smallest
convex cone $K_v$ at $v$, which contains $\Omega$. Interestingly, the latter
integrals are computed via a Fourier transform technique.

In general, when $\Omega$ has a half-space representation,
the methods described in B\"ueler et al \cite{bueler}
have a computational complexity that is exponential in $n$,
the dimension of the underlying affine space. While those methods work
well for relatively small $n$ and possibly large $m$, they become very
time-consuming and even fail for large (or even not 
so large) $n$.  This was the motivation for an alternative
method that could work in the ``dual'' context of possibly large $n$ 
and relatively small $m$.

Here we suppose given a half-space representation of $\Omega$.
The alternative method that we propose is conceptually 
very simple (as well as the computations involved) and differs from
previous existing methods.
The idea is to consider the volume of $\Omega=\{x\geq 0;Ax\leq b\}$ as
a function  $g: \R^m\to\R_+$ of the right-hand side $b\in\R^m$ for which
we provide a simple {\it explicit} expression of its
{\it Laplace transform} $G:\C^m\to\C$ in closed form. It then suffices
to apply the {\it inverse} Laplace transform to $G$, which, in the present context,
can be done efficiently by repeated applications of 
Cauchy 's Residue Theorem for the evaluation of one-dimensional complex integrals.
We propose and describe two such algorithms.

As already mentioned, the $O(n^m)$ computational complexity of both
algorithms makes the method especially attractive for large
$n$ and relatively small $m$, when the other methods
with computational complexity $O(m^n)$ would fail.
This method can also be viewed as ``dual'' of the latter methods which work in the original
space $\R^n$ with the matrix $A$, as we instead work in the space
$\R^m$ of ``dual'' variables associated with the constraints, and the cone
$\{u\geq 0,\,A'u\geq 0\}$ (via the Laplace transform), which explains
the computational complexity $O(n^m)$ (in lieu of $O(m^n)$).

\section{Main result}

Let $e_i:=(1,1,\cdots)$ be the unit vector of $\R^i$ for $i\geq 1$.
Let $y\in\R^m$ and $A\in\R^{m\times n}$ be real-valued matrices such
that the convex polyhedron 
\begin{equation}
\label{poly}
\Omega(y)\,:=\,\{x\in\R^n_+\,|\,Ax\,\leq\,y\}
\end{equation}
is compact, that is, $\Omega(y)$ is a convex polytope. The notation
$\R_+$ stands for the semi-closed interval $[0,\infty)\subset\R$.

Now consider the function $g:\R^m\to\,\R$ defined by
\begin{equation}
\label{b1}
y\,\mapsto \,g(y)\,:=\,\int_{\Omega(y)}\,dx \,=\,\mbox{vol}(\Omega(y)),
\end{equation}
and let $G:\C^m\,\to\,\C$ be its $m$-dimensional two-sided Laplace 
transform, that is,
\begin{equation}
\label{b2}
\lambda \,\mapsto\,G(\lambda)\,:=\,\int_{\R^m}\e ^{-\langle \lambda,
y\rangle}g(y)\,dy.
\end{equation}

We have the following result :
\begin{theorem}
\label{th1}
Let $\Omega(y)$ be the convex polytope in
(\ref{poly}), functions $g$ and $G$ are defined as in (\ref{b1}) and
(\ref{b2}) respectively, and assume that $x=0$ is the only solution 
of the system $\{x\geq 0,\;Ax\leq 0\}$. Then :
\begin{equation}
\label{b3}
G(\lambda)\,=\,\frac{1}{\prod_{i=1}^m\lambda_i}\times
\frac{1}{\prod_{j=1}^n(A'\lambda)_j},\hspace{1cm}
\begin{array}{c}
\Re(\lambda)\,>\,0,\\
\Re(A '\lambda)\,>\,0.
\end{array}
\end{equation}
Moreover,
\begin{equation}
\label{b4}
g(y)\,=\,\frac{1}{(2\pi i)^m}\int_{c_1-i\infty}^{c_1+i\infty}\cdots
\int_{c_m-i\infty}^{c_m+i\infty}\e ^{\langle \lambda ,y\rangle}
G(\lambda)\,d\lambda
\end{equation}
where the real constants $c>0$ satisfies $A'c>0$.
\end{theorem}
\begin{proof}
Apply the definition (\ref{b2}) of $G$, to obtain :
\begin{eqnarray*}
G(\lambda)&=&\int_{\R^m}\e ^{-\langle \lambda,y\rangle}
\left[\int_{x\geq 0,\;Ax\leq y}\,dx\right]\,dy\\
&=&\int_{\R^n_+}\left[\int_{y\geq Ax}\e^{-\langle
\lambda,y\rangle}\,dy\right]\,dx\\
&=&\frac{1}{\prod_{i=1}^m\lambda_i}
\int_{\R^n_+}\e ^{-\langle A
'\lambda ,x\rangle}\,dx,\hspace{0.5cm}\Re(\lambda)\,>\,0\\
&=&\frac{1}{\prod_{i=1}^m\lambda_i}\times \frac{1}{\prod_{j=1}^n(A
'\lambda)_j},
\hspace{0.5cm}\mbox{ with }
\left\{\begin{array}{rcl}
\Re(\lambda)&>&0\\
\Re(A '\lambda )&>&0
\end{array}\right. .
\end{eqnarray*}
And (\ref{b4}) is obtained by a direct application of the 
inverse Laplace transform. It remains to show that, indeed, the domain
$\{\Re(\lambda)>0,\,\Re(A'\lambda)>0\}$ is nonempty. However, this fact
follows from a special version of Farkas' lemma due to Carver (see e.g.
Schrijver \cite[(33), p. 95]{schrijver}), which (adapted to the 
present context) states that $\{u>0,\,A'u>0\}$ has an admissible solution 
$u\in\R^m$ if and only if $(x,y)=0$ is the only solution of the system
$\{Ax+y=0,\,x\geq 0,\, y\geq 0\}$. In other words, $x=0$ is the 
only solution of $\{x\geq 0,\,Ax\leq 0\}$.
\end{proof}
\begin{remark}
\label{remark1}
A necessary and sufficient condition for $\Omega(y)$ to be compact
is that there exists some $u\in\R^m_+$ such that $A 'u\geq e_n$.
This is a consequence of the well-known Farkas Lemma. 
\end{remark}

As mentioned above, computing $g(y)$ via (\ref{b4}) reduces
to computing the Laplace inverse of $G(\lambda)$. In our case, 
this can be done quite efficiently even for large size problems.
We first slightly modify our problem as follows :

Suppose that we want to compute the volume of the convex polytope
$\{x\geq 0;\:Ax\leq b\}$ with $b>0$, that is, we want to evaluate $g(y)$ 
at the point $y:=b\in\R^m$. We may and shall assume, without loss of 
generality, that $y_i=1$ for every $i=1,\dots m$. Otherwise, just divide 
by $y_i>0$.

The problem is then to compute the value $h(1)$ of the function
$h:\R_+\to\R$ given by
\begin{equation}
\label{b44}
h(z)=g(e_mz)=\frac{1}{(2\pi i)^m}\int_{c_1-i\infty}^{c_1+i\infty}
\cdots\int_{c_m-i\infty}^{c_m+i\infty}\e^{z\langle\lambda,e_m\rangle}\,
G(\lambda)\,d\lambda ,
\end{equation}
where the real vector $0<c\in\R^m$ satisfies $A'c>0$. Computing the
complex integral (\ref{b44}) can be done in two ways that are explored below. 
We do it directly in  \S \ref{direct}
by integrating first with respect to (w.r.t.) $\lambda_1$, then w.r.t.
$\lambda_2$, etc..., or indirectly in \S \ref{associated},
by first computing the one-dimensional Laplace 
transform $H$ of $h$ and then computing the Laplace inverse of $H$.

\section{The direct method}
\label{direct}

To better understand the methodology behind the direct method
and for illustration purpose, consider the case of a
convex polytope $\Omega$ with only two ($m=2$) 
nontrivial constraints. 

\subsection{The $m=2$ non trivial constraints example}
Let $A\in \R^{2\times n}$ be such
that $x=0$ is the only solution of $\{x\geq0,\,Ax\leq0\}$. Moreover, 
suppose that $A ':=[a\,|\,b]$ with $a,b\in\R^n$. 
For ease of exposition, assume that 

\begin{itemize}
\item $a_jb_j\neq 0$ and $a_j\neq b_j$ for all $j=1,2,\dots n$.
\item $a_j/b_j\neq a_k/b_k$ for all $j,k=1,2,\dots n$
\end{itemize}
Observe that these assumptions are satisfied with probability one if we add to 
every coefficient $a_i,b_i$ a perturbation $\epsilon$, randomly
generated under a uniform distribution on $[0,\bar{\epsilon}]$, 
with $\bar{\epsilon}$ very small.

Then :
\[G(\lambda)\,=\,\frac{1}{\lambda_1\lambda_2}\times
\frac{1}{\prod_{j=1}^n (a_j\lambda_1+b_j\lambda_2)},
\hspace{1cm}\left\{ \begin{array}{rcl}\Re(\lambda)&>&0\\
\Re(a\lambda_1+b\lambda_2)&>&0\end{array}\right. .\]

Next, fix $c_1$ and $c_2>0$ such that 
$a_jc_1+b_jc_2>0$ for every $j=1,2,\dots n$, and compute
the integral (\ref{b44}) as follows. 
We first evaluate the integral
\begin{equation}
\label{I_1}
I_1\,=\,\frac{1}{2\pi i}\int_{c_1-i\infty}^{c_1+i\infty}\frac{\e^{z\lambda_1}}
{\lambda_1\prod_{j=1}^n (a_j\lambda_1+\lambda_2b_j)}\,d\lambda_1,
\end{equation}
using classical Cauchy 's residue technique. That is, we: (a) close the path
of integration adding a semicircle $\Gamma$ of radius $R$ large enough, (b) 
evaluate the closed integral using Cauchy 's Residue Theorem
\cite[Theor. 2.2, p. 112]{conway}, and (c) show 
that the integral along $\Gamma$ converges to zero when $R\rightarrow\infty$.

Now, since we are integrating w.r.t. $\lambda_1$ and we want to
evaluate $h(z)$ at $z=1$, we must add the semicircle $\Gamma$ on the 
left side of the integration path $\Re(\lambda_1)=c_1$ because 
$\e^{y\lambda_1}$ converges to zero when $\lambda_1\rightarrow-\infty$. 
Therefore, we must consider {\it only} the poles of 
$G(\lambda_1,\cdot)$ whose real part is strictly less than $c_1$
(with $\lambda_2$ being fixed). Then, the evaluation of (\ref{I_1})
follows easily, and
\[I_1\,=\,\frac{1}{\lambda_2^n\,\prod_{j=1}^nb_j}
+\sum_{j=1}^n\frac{-e^{-(b_j/a_j)z\lambda_2}}
{b_j\lambda_2^n\,\prod_{k\neq j}(-a_kb_j/a_j+b_k)}.\]
Recall that $\Re(-\lambda_2b_j/a_j)=-c_2b_j/a_j<c_1$ 
for each $j=1,2,\dots n$, and $G(\lambda_1,\cdot)$ 
has only poles of first order (with $\lambda_2$ being fixed).

Therefore, 
\begin{eqnarray*}
h(z)&=&\frac{1}{2\pi i}\int_{c_2-i\infty}^{c_2+i\infty}
\frac{\e^{z\lambda_2}}{\lambda_2}\;I_1\;d\lambda_2\\
&=&\frac{1}{2\pi i}\int_{c_2-i\infty}^{c_2+i\infty}\frac{\e^{z\lambda_2}}
{\lambda_2^{n+1}\,\prod_{j=1}^nb_j}\;d\lambda_2\;-\\
&&-\;\sum_{j=1}^m\frac{1}{2\pi i}\int_{c_2-i\infty}^{c_2+i\infty}
\frac{a_j^n\;\e^{(1-b_j/a_j)z\lambda_2}}
{\lambda_2^{n+1}a_jb_j\,\prod_{k\neq j}(b_ka_j-a_kb_j)}\;d\lambda_2.\\
\end{eqnarray*}

These integrals must be evaluated according to whether $(1-b_j/a_j)y$ is 
positive or negative. Thus, recalling that $z>0$, each integral is equal to 

\indent
- its  residue at the pole $\lambda_2=0<c_2$ when $1-b_j/a_j$ is
positive, and 

\indent
- zero if $1-b_j/a_j$ is negative because there is no pole on 
the right side of $\Re(\lambda_2)=c_2$. 

That is,
\begin{equation}
\label{g(y)}
h(z)\,=\,\frac{z^n}{n\mbox{!}}\left[\frac{1}{\prod_{j=1}^nb_j}
-\sum_{b_j/a_j<1}\frac{(a_j-b_j)^n}{a_jb_j\prod_{k\neq j}
(b_ka_j-a_kb_j)}\right]\;.
\end{equation}
Observe that the formula is not symmetrical in the parameters $a,b$. 
This is because we have chosen to integrate first w.r.t.
$\lambda_1$; and the set $\{j\,|\,b_j/a_j<1\}$ is different from 
$\{j\,|\,a_j/b_j>1\}$, which would have been considered had we 
integrated first w.r.t. $\lambda_2$. In the latter case, 
we would have obtained
\begin{equation}
\label{secondform}
h(z)\,=\,\frac{z^n}{n\mbox{!}}\left[\frac{1}{\prod_{j=1}^na_j}
-\sum_{a_j/b_j<1}\frac{(b_j-a_j)^n}{a_jb_j\prod_{k\neq j}
(a_kb_j-b_ka_j)}\right],
\end{equation}
which is (\ref{g(y)}) by interchanging $a$ and $b$. Moreover, moving terms
around, we get for free the following identity 
\begin{equation}
\label{ide}
\sum_{j=1}^n\frac{(a_j-b_j)^n}{a_jb_j\prod_{k\neq j}(b_ka_j-a_kb_j)}
\;=\;\frac{1}{\prod_{j=1}^nb_j}\;-\;\frac{1}{\prod_{j=1}^na_j}\;.
\end{equation}

\subsection{The direct method algorithm}
\label{directsub}

The above methodology easily extends to an arbitrary number $m$ of
non trivial constraints.  One evaluates the
integral of the right-hand side of (\ref{b44}) by integrating first
w.r.t. $\lambda_1$, then w.r.t. $\lambda_2$, and so on.
The resulting algorithm can be described with a tree
of depth $m+1$ ($m+1$ ``levels'').  Let $0 <c\in\R^m$ be such that $A
'c>0$.

- Level $0$ is the root of the tree.

- Level $1$ is the integration w.r.t. $\lambda_1$ and consists of 
at most $n+1$ nodes associated with the poles $\lambda_1:=\rho_j^1$, $j=1,\dots n+1$, of the
rational function $\prod_i\lambda_i ^{-1}\prod_j (A '\lambda)^{-1}_j$, seen as a
function of $\lambda_1$ only. 
By the assumption on $c$, there is no pole $\rho_j^1$ on the line $\Re(\lambda_1)=c_1$.
By Cauchy 's Residue Theorem, only the poles at the left side of 
the integration path $\Re(\lambda_1)=c_1$, say $\rho_j^1$, $j\in I_1$,
are selected.

- Level $2$ is the integration w.r.t. $\lambda_2$. After integration
w.r.t. $\lambda_1$, {\it each} of the poles $\rho_j^1$, $j\in I_1$,
generates a rational function of $\lambda_2,\lambda_3,\dots,\lambda_m$,
which, seen as a function of $\lambda_2$ only, has at most
$n+1$ poles $\rho_i^2(j)$, $i=1,\dots n+1$, $j\in I_1$. Thus level $2$
has at most $(n+1)^2$ nodes associated with the poles $\rho_i^2(j)$.
Assuming no pole $\rho_i^2(j)$ on the line $\Re(\lambda_2)=c_2$,
by Cauchy 's Residue Theorem, only the poles 
$\rho_i^2(j)$, $(j,i)\in I_2$, located on the correct side of the
integration path $\Re(\lambda_2)=c_2$ are selected.

- Level $k$, $k\leq m$, consists of
at most $(n+1)^k$ nodes associated with the poles 
$\rho_s^k(i_1,i_2,\dots i_{k-1})$, $(i_1,i_2,\dots i_{k-1})\in I_{k-1}$,
$s=1,\dots n+1$, of some rational functions of $\lambda_k,\dots ,\lambda_m$,
seen as functions of $\lambda_k$ only. Assuming no pole
on the line $\Re(\lambda_k)=c_k$, only the poles 
$\rho_{i_k}^k(i_1,i_2,\dots i_{k-1})$, $(i_1,i_2,\dots i_k)\in I_{k}$,
located on the correct side of the integration path
$\Re(\lambda_k)=c_k$, are selected. And so on.

The last level $m$ consists of at most $(n+1)^m$ nodes and the
integration w.r.t. $\lambda_m$ is trivial as it amounts to evaluate
integrals of the form 
\[(2\pi i)^{-1}\int_{c_m-i\infty}^{c_m+i\infty} A\lambda_m^{-(n+1)}\e
^{\alpha z\lambda_m}d\lambda_m,\]
for some coefficients $A$, $\alpha$, which yields $A(\alpha
z)^n/n\mbox{!}$ for those $\alpha >0$. Summing up over all
the nodes provides the desired value.

Only simple elementary arithmetic operations are needed to
compute the nodes at each level, as in Gauss elimination for solving
linear systems. Therefore, the computational complexity is easily seen to
be $O(n^m)$. 

However, some care must be taken
with the choice of the integration paths as we assume
that at each level $k$ there is no pole on the integration path $\Re(\lambda_k)=c_k$. This
issue is discussed in \S \ref{path}. The algorithm is illustrated on the following simple example
with $n=2, m=3$.
\vspace{0.2cm}

\noindent
{\bf Example:} Let $\Omega(ze_2)\subset\R^2$ be the polytope
\[\Omega(ze_2)\,:=\,\{x\in\R^2_+\,|\,x_1+x_2\,\leq\,z; -2x_1+2x_2\leq
z\,;\:2x_1-x_2\leq z\},\]
whose area is $17z^2/48$.

Choose $c_1=3$, $c_2=2$ and $c_3=1$, 
so that $c_1>2c_2-2c_3$ and $c_1>c_3-2c_2$.
\[h(z)\,=\,\frac{1}{(2\pi i)^3}\int_{c_1-i\infty}^{c_1+i\infty}\dots
\int_{c_3-i\infty}^{c_3+i\infty}\e^{(\lambda_1+\lambda_2+\lambda_3)z}
G(\lambda)\,d\lambda,\]
with
\[G(\lambda)\,=\,\frac{1}{\lambda_1\lambda_2\lambda_3
(\lambda_1-2\lambda_2+2\lambda_3)(\lambda_1+2\lambda_2-\lambda_3)}.\]

Integrate first w.r.t. $\lambda_1$; that is, evaluate the 
residues at the poles $\lambda_1=0$, $\lambda_1=2\lambda_2-2\lambda_3$ 
and $\lambda_1=\lambda_3-2\lambda_2$ because $0<z$, $0<c_1$, 
$\Re(2\lambda_2-2\lambda_3)<c_1$ and
$\Re(\lambda_3-2\lambda_2)<c_1$. We obtain
\[h(z)=\frac{1}{(2\pi i)^{2}}\int_{c_2-i\infty}^{c_2+i\infty}
\int_{c_3-i\infty}^{c_3+i\infty}I_2+I_3+I_4\;d\lambda_2\,d\lambda_3,\]
where
\begin{eqnarray*}
I_2&=&\frac{-\e ^{(\lambda_2+\lambda_3)z}}
{2\lambda_2\lambda_3(\lambda_3-\lambda_2)(\lambda_3-2\lambda_2)},\\
I_3&=&\frac{\e^{(3\lambda_2-\lambda_3)z}}
{6\lambda_2\lambda_3(\lambda_3-\lambda_2)(\lambda_3-4\lambda_2/3)},\\
I_4&=&\frac{\e ^{(2\lambda_3-\lambda_2)z}}
{3\lambda_2\lambda_3(\lambda_3-2\lambda_2)(\lambda_3-4\lambda_2/3)}.
\end{eqnarray*}
Next, integrate $I_2$ w.r.t. $\lambda_3$. We must consider 
the poles on the left side of $\Re(\lambda_3)=1$, that is, the pole 
$\lambda_3=0$ because $\Re(\lambda_2)=2$. Thus, we get 
$-\e^{z\lambda_2}/4\lambda_2^3$, and the next integration w.r.t.
$\lambda_2$ yields $-z^2/8$.

When integrating $I_3$ w.r.t. $\lambda_3$,
we have to consider the poles $\lambda_3=\lambda_2$ and
$\lambda_3=4\lambda_2/3$, on the right side of $\Re(\lambda_3)=1$;
and we get
\[\frac{-1}{\lambda_2^3}\left[-\frac{\e^{2z\lambda_2}}{2}+
\frac{3\e^{z\lambda_25/3}}{8}\right]\,.\]
Recall that the path of integration has a negative orientation, so we 
have to consider the negative value of residues. The next integration 
w.r.t. $\lambda_2$ yields $z^2(1-25/48)$.

Finally, when integrating $I_4$ w.r.t. $\lambda_3$,
we must consider only the pole $\lambda_3=0$, 
and we get $\e^{-z\lambda_2}/8\lambda_2^3$;
the next integration w.r.t. $\lambda_2$ yields zero. Hence,
adding up the above three partial results, yields
\[h(z)\,=\,z^2\left[\frac{-1}{8}+1-\frac{25}{48}\right]\,
=\,\frac{17\,z^2}{48},\]
which is the desired result.
\vspace{0.2cm}

\subsection{Paths of integration}
\label{path}

In choosing the integration paths $\Re(\lambda_k)=c_k$,
$k=1,\dots m$, we must determine
a vector $0<c\in\R^m$ such that $A'c>0$. However, this may not be enough when 
we want to evaluate the integral (\ref{b44}) by repeated
applications of Cauchy's Residue Theorem. Indeed, we have seen in the
tree description of the algorithm (cf. \S \ref{directsub}), that at
each level $k>1$ of the tree (integration
w.r.t $\lambda_k$), one {\it assumes } that there is {\it no} pole on
the integration path $\Re(\lambda_k)=c_k$.

For instance, had we set $c_1=c_2=c_3=1$ (instead of $c_1=3$, $c_2=2$ 
and $c_1=1$) in the above example, we could not use 
Cauchy's Residue Theorem to integrate $I_2$ or $I_3$ because we would 
have the pole $\lambda_2=\lambda_3$ exactly on the path of integration
(recall that $\Re(\lambda_2)=\Re(\lambda_3)=1$); fortunately, this
case is {\it pathological} as it happens with probability zero in a set of problems
with randomly generated data $A\in \R^{m\times n}$ and, therefore,
this issue could be ignored in practice. However, for the sake of mathematical rigor,
in addition to the constraints $c>0$ and $A 'c>0$,
the vector $c\in\R^m$  must satisfy additional constraints to avoid the
above mentioned pathological problem. We next describe one way to
proceed to ensure that $c$ satisfies these additional constraints.

In \S \ref{directsub} we have described the algorithm as a tree of depth $m$
(level $i$ being the integration w.r.t. $\lambda_i$) where
each node has at most $n+1$ descendants (one descendant for each pole on the
correct side of the integration path $\Re(\lambda_i)=c_i$). 
The volume is then the summation of all partial results obtained at
each leaf of the tree (that is, each node of level $m$). We next describe how
to ``perturbate'' {\it on-line} the initial vector $c\in\R^m$ if at
some level $k$ of the algorithm
there is a pole on the corresponding integration path $\Re(\lambda_k)=c_k$.

{\it - Step 1. Integration w.r.t. $\lambda_1$.}
Choose a real vector $c:=(c_1^1,\cdots,c_m^1)>0$ such that
$A'c>0$ and integrate (\ref{b44}) along 
the line $\Re(\lambda_1)=c_1^1$. From Cauchy 's Residue Theorem,
this is done by selecting the (at most $n+1$) poles
$\lambda_1:=\rho_j^1$, $j\in I_1$, located on the left-side of the integration path $\Re(\lambda_1)=c_1^1$.
Each pole $\rho_j^1, j=1,\dots n+1$ (with $\rho_j^1:=0$) is a linear combination 
$\beta_{j2}^{(1)}\lambda_2+\ldots+\beta_{jm}^{(1)}\lambda_m$ with real coefficients
$\{\beta_{jk}^{(1)}\}$, because $A$ is a real-valued matrix. 
Observe that by the initial choice of $c$,
\[\delta_1\,:=\,\min_{j=1,\dots n+1}\,\vert
c_1^1-\sum_{k=2}^m\beta_{jk}^{(1)}c_k^1\vert\,>\,0.\]

{\it - Step 2. Integration w.r.t. $\lambda_2$}. 
For each of the poles $\rho_j^1$, $j\in I_1$, selected at step 1, and
after integration w.r.t. $\lambda_1$, we now have to consider a 
rational function of $\lambda_2$ with at most $n+1$ poles
$\lambda_2:=\rho_i^2(j):=\sum_{k=3}^m\beta_{ik}^{(2)}(j)\lambda_k$, $i=1,\dots
n+1$. If
\[\delta_2\,:=\,\min_{j\in I_1}\min_{i=1,\dots n+1}\,\vert
c_2^1-\sum_{k=3}^m\beta_{ik}^{(2)}(j)c_k^1\vert\,>\,0,\]
then integrate w.r.t. $\lambda_2$ on the line
$\Re(\lambda_2)=c_2^1$. Otherwise, if $\delta_2=0$
we  set $c_2^2:=c_2^1+\epsilon_2$ and
$c_k^2:=c_k^1$ for all $k\neq 2$, 
by choosing $\epsilon_2>0$ small enough to ensure
that 
\[\begin{array}{lrcl}
\mbox{(a)}&A 'c^2 &>&0\\
\mbox{(b)}&\delta_2\,:=\,\displaystyle{\min_{j\in I_1}\min_{i=1,\dots n+1}\vert
c_2^2-\sum_{k=3}^m\beta_{ik}^{(2)}(j)c_k^2\vert}&>&0\\
\mbox{(c)}&\displaystyle{\max_{j=1,\dots n+1}\vert \beta_{j2}^{(1)}\epsilon_2\vert}&<&\delta_1
\end{array}\]
The condition (a) is basic whereas (b) ensures that 
there is no pole on the integration path $\Re(\lambda_2)=c_2^2$.
Moreover, what has been done in step $1$ remains valid because from (c),
$c_1^2-\sum_{k=2}^m\beta_{jk}^{(1)}c_k^2$ has the same sign
as $c_1^1-\sum_{k=2}^m\beta_{jk}^{(1)}c_k^1$, and, therefore, none of
the poles $\rho_j^1$ has crossed the integration path
$\Re(\lambda_1)=c_1^1=c_1^2$, that is, the set $I_1$ is unchanged.

Then integrate w.r.t. $\lambda_2$ on the line $\Re(\lambda_2)=c_2^2$,
which is done via Cauchy 's Residue Theorem by
selecting the (at most $(n+1)^2$) poles $\rho^2_i(j)$,
$(j,i)\in I_2$, located at the left or the right of the line
$\Re(\lambda_2)=c_2^2$,  depending on the sign of the coefficient of
the argument in the exponential.

{\it - Step 3. Integration w.r.t. $\lambda_3$}. 
Likewise, for each of the poles $\rho_i^2(j)$, $(j,i)\in I_2$, selected
at step 2, we now have to consider a rational function of $\lambda_3$ with at most $n+1$ poles
$\rho_s^3(j,i):=\sum_{k=4}^m\beta_{sk}^{(3)}(j,i)\lambda_k$, $s=1,\dots n+1$. If
\[\delta_3\,:=\,\min_{(j,i)\in I_2}\min_{s=1,\dots n+1}\,\vert
c_3^2-\sum_{k=4}^m\beta_{sk}^{(3)}(j,i)c_k^2\vert\,>\,0,\]
then integrate w.r.t. $\lambda_3$ on the line $\Re(\lambda_3)=c_3^2$. 
Otherwise, if $\delta_3 =0$, set
$c_3^3:=c_3^2+\epsilon_3$ and
$c_k^3:=c_k^2$ for all $k\neq 3$, by choosing $\epsilon_3>0$
small enough to ensure that 
\[\begin{array}{lrcl}
\mbox{(a)}&A 'c^3 &>&0\\
\mbox{(b)}&
\delta_3\,:=\,\displaystyle{\min_{(j,i)\in I_2}\min_{s=1,\dots n+1}}
\vert c_3^3-\sum_{k=4}^m\beta_{sk}^{(3)}(j,i)c_k^3\vert&>&0\\
\mbox{(c)}&\displaystyle{\max_{j\in I_1}\max_{i=1,\dots n+1}}\vert
\beta_{i3}^{(2)}(j)\epsilon_3\vert&<&\delta_2\\
\mbox{(d)}&
\displaystyle{\max_{j=1,\dots n+1}}\vert\beta_{j2}^{(1)}\epsilon_2+\beta_{j3}^{(1)}\epsilon_3\vert&<&\delta_1
\end{array}\]
As in previous steps, condition (a) is basic.
The condition (b) ensures that there is no pole on the integration path
$\Re(\lambda_3)=c_3^3$. Condition (c) (resp. (d)) ensures that 
none of the poles $\rho ^2_i(j)$ considered at step $2$ (resp. none of the poles $\rho
^1_j$ considered at step $1$) has crossed 
the line $\Re(\lambda_2)=c_2^3=c_2^2$ (resp. the line
$\Re(\lambda_1)=c_1^3=c_1^1$). That is, both sets $I_1$ and $I_2$ are unchanged.

Then integrate w.r.t. $\lambda_3$ on
the line $\Re(\lambda_3)=c_3^3$, which is done by selecting the
(at most $(n+1)^3$) poles $\rho_s(j,i)$, $(j,i,s)\in I_3$, located at the left or right of the line
$\Re(\lambda_3)=c_3^3$, depending on the sign of the argument in the exponential.

And so on. It is important to notice that the $\epsilon_k$'s and
$c_k^k$'s play no (numerical) role in the integration itself. They
are only used to (i) ensure
the absence of a pole on the integration path $\Re(\lambda_k)=c_k^k$,
and (ii) to locate the poles on the left or the right of the
integration path. Their numerical value (which can be very small) has
no influence on the computation.

\section{The associated transform algorithm}
\label{associated}

An alternative to the direct method
permits to avoid evaluating integrals of exponential functions in 
(\ref{b44}) by making the following simple change of variable. Let
$\lambda_1=p-\sum_{j=2}^m\lambda_j$ and $d=\sum_{j=1}^mc_j$ in
(\ref{b44}), so that
\[h(z)=\frac{1}{(2\pi i)^m}\int_{c_m-i\infty}^{c_m+i\infty}\dots
\int_{c_2-i\infty}^{c_2+i\infty}\left[\int_{d-i\infty}^{d+i\infty}
\e^{zp}\widehat{G}\;dp\right]\;d\lambda_2\dots d\lambda_m,\]
where
\begin{equation}
\label{aso}
\widehat{G}\;=\;G(p-\sum_{j=2}^m\lambda_j,\lambda_2,\dots,\lambda_m).
\end{equation}

We can rewrite $h(z)$ as
\begin{eqnarray}
\label{b55}
h(z)&=&\frac{1}{2\pi i}\int_{d-i\infty}^{d+i\infty}\;\e^{zp}H(p)dp,
\hspace{1cm}\hbox{with}\\
\label{b5}
H(p)&:=&\frac{1}{(2\pi i)^{m-1}}\int_{c_2-i\infty}^{c_2+i\infty}\dots
\int_{c_m-i\infty}^{c_m+i\infty}\widehat{G}\;d\lambda_2\dots d\lambda_m.
\end{eqnarray}

Recall that $G(\lambda)$ is well defined on the domain 
$\Re(\lambda)>0$ and $\Re(A'\lambda)>0$; moreover, the real vector $c$
is taken in this domain. Hence, the domain of definition 
of $H(p)$ is given by the condition 
\[(\Re(p)-\sum_{j=2}^mc_j,c_2,\dots,c_m)\;
\in\;\{y\in\R^m\,|\,y\,>\,0,\,A'y\,>\,0\}.\]

On other hand, recall that the system $\{x\geq0,\,Ax\leq0\}$ has only one
solution $x=0$ (see the hypotheses of Theorem \ref{th1}). Hence, the function 
$h(z)$ is identically equal zero when $z\leq 0$ (see (\ref{b1}) 
and (\ref{b44})). 
Therefore, $H(p)$ is the one-sided Laplace transform of $h(z)$.
Moreover, it is also easy to see that there exists 
a real constant $C$ such that $h(z)=z^nC/n!$ when $z\geq0$. Therefore,
\[H(p)=C/p^{n+1}\]
and the main problem completely reduces to evaluating
the constant $C=h(1)n!$ by integrating $\widehat{G}$ in (\ref{b5}).

Notice that we only need to evaluate $m-1$ integrals.
The function $H(p)$ is called the {\it associated transform} of $G(\lambda)$.

Again, the integral (\ref{b5}) can be computed via repeated
applications of Cauchy 's Residue  Theorem
(and as in the direct method algorithm of \S \ref{direct},
some care is needed with the domain of integration 
and the location of the poles of $\widehat{G}$).
The method is illustrated on the same example of two non trivial 
constraints  ($m=2$) already considered at the beginning of \S \ref{direct}.

\subsection{The $m=2$ non trivial constraints example}
Let $A\in \R^{2\times n}$ such that $x=0$ is the only solution of 
$\{x\geq0,\,Ax\leq0\}$. Write $A ':=[a\,|\,b]$ 
with $a,b\in\R^n$. To compare with the direct method,
and as in the beginning of \S \ref{direct},
assume that $a_jb_j\neq 0$ for all $j=1,\dots$ and 
$a_j/b_j\neq a_k/b_k$ for all $j\neq k$.

Then :
\[G(\lambda)\,=\,\frac{1}{\lambda_1\lambda_2}\times
\frac{1}{\prod_{j=1}^n (a_j\lambda_1+b_j\lambda_2)},
\hspace{1cm}\begin{array}{c}\Re(\lambda)>0,\\
\Re(A'\lambda)>0.\end{array}\]
Fix $\lambda_2=p-\lambda_1$ and choose a real constant $c_1>0$ such that 
the system of inequalities $\Re(p)>c_1$ and $(a_j-b_j)c_1+b_j\Re(p)>0$ 
for all $j=1,\dots n$ has a solution. We already know that there is at least 
one vector $u=(c_1,\Re(p)-c_1)$ such that $u>0$ and 
$A'u>0$. We obtain $H(p)$ by integrating $G(\lambda_1,p-\lambda_1)$ 
w.r.t. $\lambda_1$, which yields

\[H(p)\,=\,\frac{1}{2\pi i}\int_{c_1-i\infty}^{c_1+i\infty}
\frac{1}{\lambda_1(p-\lambda_1)}\times
\frac{1}{\prod_{j=1}^n ((a_j-b_j)\lambda_1+b_jp)}\,d\lambda_1,\]

Next, we need to determine which poles of $G(\lambda_1,p-\lambda_1)$ are
on the left (right) side of the integration path $\Re(\lambda_1)=c_1$
in order to apply Cauchy's Residue theorem. Let $J_+=\{j|a_j>b_j\}$, 
$J_0=\{j|a_j=b_j\}$ and $J_-=\{j|a_j<b_j\}$. Then, the poles on the left
side of $\Re(\lambda_1)=c_1$ are $\lambda_1=0$ and 
$\lambda_1=-b_jp/(a_j-b_j)$ for all $j\in J_+$ because 
$-b_j\Re(p)/(a_j-b_j)<c_j$. Besides, the poles on the right side of 
$\Re(\lambda_1)=c_1$ are $\lambda_1=p$ and 
$\lambda_1=-b_jp/(a_j-b_j)$ for all $j\in J_-$. Finally, notice that
$G(\lambda_1,p-\lambda_1)$ has only poles or first order. 

Hence, computing the residues of poles on the left side of
$\Re(\lambda_1)=c_1$, yields
\begin{eqnarray*}
H(p)&=&\frac{1}{\prod_{j\in J_0}pb_j}\,\left[
\frac{1}{p\prod_{j\not\in J_0}pb_j}\;+\right.\\
&&\left.+\;\sum_{j\in J_+}\frac{-(a_j-b_j)^{n-|J_0|}}
{p^2\,a_jb_j\,\prod_{k\not\in J_0, k\neq j}(-pb_ja_k+pa_jb_k)}\right]\\
&=&\frac{1}{p^{n+1}}\left[\frac{1}{\prod_{j=1}^nb_j}\;-\;
\sum_{a_j/b_j>1}\frac{(a_j-b_j)^n}{a_jb_j\,\prod_{k\neq j}
(a_jb_k-b_ja_k)}\right],
\end{eqnarray*}
and one retrieves (\ref{g(y)}) when we take $J_0$ to be an empty 
set, in other words, when its cardinality $|J_0|=0$. Now, computing
the negative value of residues of poles on the right side of 
$\Re(\lambda_1)=c_1$ (we need to take the negative value because the
path of integration has a negative orientation), yields
\[H(p)=\frac{1}{p^{n+1}}\left[\frac{1}{\prod_{j=1}^na_j}\;-\;
\sum_{b_j/a_j>1}\frac{(b_j-a_j)^n}{a_jb_j\,\prod_{k\neq j}
(b_ja_k-a_jb_k)}\right],\]
and we also retrieve \ref{secondform}.

\subsection{The associated transform algorithm}

As for the direct method algorithm, the above methodology easily
extends to an arbitrary number $m$ of nontrivial constraints. 
The algorithm also consists of $m$ 
(one-dimensional integration) steps. At each step,
the several one-dimensional complex integrals are evaluated by application of 
Cauchy 's Residue Theorem \cite[Theor. 2.2, p. 112]{conway}.
For same reasons as in the direct method,
the computational complexity is easily seen to be $O(n^m)$.

The general case is better illustrated on the same example as in 
\S \ref{directsub}. Again, to avoid the case of poles on the
integration path in pathological examples, some care is 
needed when one specifies the integration path at each step of the algorithm.

Let $\Omega(ze_2)\subset\R^2$ be the polytope
\[\Omega(ze_2)\,:=\,\{x\in\R^2_+\,|\,x_1+x_2\,\leq\,z; -2x_1+2x_2\leq
z\,;\:2x_1-x_2\leq z\},\]
whose area is $17z^2/48$.

We can choose $\lambda_3=p-\lambda_2-\lambda_1$ and $c_1=c_2=1$ such 
that $\Re(p)>2$, $2\Re(p)>5$ and $\Re(p)<5$; and so
\[H(p)\,=\,\frac{1}{(2\pi i)^2}\int_{1-i\infty}^{1+i\infty}
\int_{1-i\infty}^{1+i\infty}M(\lambda ,p)\,d\lambda_1\,d\lambda_2,\]
with
\[M(\lambda,p)\,=\,\frac{1}{\lambda_1\lambda_2\,(p-\lambda_1-\lambda_2)
(2p-\lambda_1-4\lambda_2)(2\lambda_1+3\lambda_2-p)}.\]

We first integrate w.r.t. $\lambda_1$. Only the real parts of the 
poles $\lambda_1=0$ and $\lambda_1=(p-3\lambda_2)/2$ are less than $1$.
Therefore, the residue of the  $0$-pole yields:
\begin{equation}
\label{c1}
\frac{1}{2\pi i}\int_{1-i\infty}^{1+i\infty}
\frac{1}{\lambda_2\,(p-\lambda_2)(2p-4\lambda_2)(3\lambda_2-p)}\,d\lambda_2,
\end{equation}
whereas the residue of the $(p-3\lambda_2)/2)$-pole yields
\begin{equation}
\label{c2}
\frac{1}{2\pi i}\int_{1-i\infty}^{1+i\infty}
\frac{4}{\lambda_2\,(p-3\lambda_2)(p+\lambda_2)(3p-5\lambda_2)}\,d\lambda_2.
\end{equation}

At this point, we have to be careful; observe that $5/2<\Re(p)<5$. However, 
we cannot put $\Re(p)=3$ because otherwise we will have a pole in the path 
of integration of (\ref{c1}) and (\ref{c2}). We thus fix
$3<\Re(p)<5$. Applying again Cauchy's Residue Theorem to (\ref{c1}) at
the pole $\lambda_2=0$ (the only one whose real part is less than one), 
yields $-1/2p^3$.

Similarly, applying Cauchy 's Residue Theorem to (\ref{c2}) at the poles
$\lambda_2=0$ and $\lambda_2=-p$ (the only ones whose real part is less 
than one), yields $4/3p^3-1/8p^3$.

We finally have that $H(p)=(4/3-1/8-1/2)/p^3=17/24p^3$, 
and so $h(z)=17z^2/48$, the desired result.

Concerning the pathological case of some poles on
the integration paths at some step of the algorithm,
the same remarks and similar remedies as for the direct method
are valid (cf. \S \ref{path}).

\section{Conclusion}

We have presented two algorithms for the 
exact computation of the volume of a convex polytope
given by its {\it half-space} representation.
The methodology behind both algorithms is conceptually simple as it
reduces to invert the Laplace
transform of the volume (considered as a function of the
right-hand-side). Both algorithms are relatively easy 
to implement (with special care for the choice of the
integration paths of the repeated one-dimensional integrals). Their $O(n^m)$ computational
complexity can make them especially attractive for large $n$ and small
$m$, when the other methods (with half-space representation of $\Omega$)
fail because of their $O(m^n)$ computational complexity.


\begin{thebibliography}{lasserre}
\bibitem{algower}
E.L. Allgower, P.M. Schmidt. Computing volumes of polyhedra,
{\it Math. Comp.} {\bf 46} (1986), 171--174.
\bibitem{barvinok}
A.I. Barvinok. Computing the volume, couting integral points and
exponentials sums, {\it Discr. Comp. Geom.} {\bf 10} (1993), 123--141.
\bibitem{bueler}
B. B\"ueler, A. Enge, K. Fukuda. Exact volume computation for
polytopes : A practical study. In:  {\it Polytopes - Combinatorics and
Computation},
G. Kalai, G. M. Ziegler, Eds., Birh\"auser Verlag, Basel, 2000.
\bibitem{cohen}
J. Cohen, T. Hickey. Two algorithms for determining volumes of
convex polyhedra, {\it J. ACM} {\bf 26} (1979), 401--414.
\bibitem{conway}
J.B. Conway. {\it Functions of a complex variable I}, 2nd ed., Springer, New
York, (1978).
\bibitem{dyer}
M.E. Dyer. The complexity of vertex enumeration methods. {\it
Math. Oper. Res.} {\bf 8} (1983), 381--402.
\bibitem{dyer2}
M.E. Dyer, A.M. Frieze. The complexity of computing the volume of a
polyhedron.
{\it SIAM J. Comput.} {\bf 17} (1988), 967--974.
\bibitem{lasserre}
J.B. Lasserre. An analytical expression and an algorithm for the
volume of a convex polyhedron in $\R^n$. {\it J. Optim. Theor. Appl.}
{\bf 39} (1983), 363--377.
\bibitem{lawrence}
J. Lawrence. Polytope volume computation, {\it Math. Comp.} {\bf 57}
(1991), 259--271.
\bibitem{schrijver}
A. Schrijver. {\it Theory of Linear and Integer Programming},
John Wiley \& Sons, Chichester, 1986.
\bibitem{von}
B. Von Hohenbalken. Finding simplicial subdivisions of polytopes,
{\it Math. Prog.} {\bf 21} (1981), 233--234.
\end{thebibliography}
\end{document}